\documentclass[reqno,centertags, 12pt]{amsart}
\usepackage{amsmath,amsthm,amscd,amssymb}
\usepackage{latexsym,verbatim}

\newcommand{\bbR}{{\mathbb{R}}}

\newcommand{\lb}{\label}

\newcommand{\beq}{\begin{equation}}
\newcommand{\eeq}{\end{equation}}
\newcommand{\ba}{\begin{align}}
\newcommand{\ea}{\end{align}}
\newcommand{\eps}{\varepsilon}
\newcommand{\del}{\delta}
\newcommand{\tht}{\theta}

\newcommand{\til}{\tilde}

\newcounter{smalllist}
\newenvironment{SL}{\begin{list}{{\rm\roman{smalllist})}}{%
\setlength{\topsep}{0mm}\setlength{\parsep}{0mm}\setlength{\itemsep}{0mm}%
\setlength{\labelwidth}{2em}\setlength{\leftmargin}{2em}\usecounter{smalllist}%
}}{\end{list}}

\DeclareMathOperator*{\dist}{dist}

\allowdisplaybreaks \numberwithin{equation}{section}

\newtheorem{theorem}{Theorem}

\newtheorem{lemma}[theorem]{Lemma}

\theoremstyle{definition}

\theoremstyle{remark}

\begin{document}
\title[Sharp Transition Between Extinction and Propagation]{Sharp Transition Between Extinction and Propagation of Reaction}

\author{Andrej Zlato\v s}

\address{ Department of Mathematics \\ University of
Wisconsin \\ Madison, WI 53706, USA \\ Email: \tt
zlatos@math.wisc.edu}


\maketitle

\begin{abstract}
We consider the reaction-diffusion equation
\[
T_t = T_{xx}  + f(T)
\]
on $\bbR$ with $T_0(x) \equiv \chi_{[-L,L]} (x)$ and
$f(0)=f(1)=0$. In 1964 Kanel' proved that if $f$ is an ignition
non-linearity, then $T\to 0$ as $t\to\infty$ when $L<L_0$, and
$T\to 1$ when $L>L_1$. We answer the open question of relation of
$L_0$ and $L_1$ by showing that $L_0=L_1$. We also determine the
large time limit of $T$ in the critical case $L=L_0$, thus
providing the phase portrait for the above PDE with respect to a
1-parameter family of initial data. Analogous results for
combustion and bistable non-linearities are proved as well.
\end{abstract}

\section{Introduction} \lb{S1}

In the present paper we consider the reaction-diffusion equation
\begin{equation} \lb{0.1}
T_t = \Delta T  + f(T)
\end{equation}
in the cylinder $\bbR\times\Omega$ where $\Omega$ is a domain in
$\bbR^{n-1}$, with Neumann boundary conditions on
$\bbR\times\partial\Omega$. The non-linear reaction term $f$ is
assumed to be Lipschitz continuous with $f(0)=f(1)=0$ and the
initial datum $T_0$ is between 0 and 1.

We will treat the case when $T_0$ is independent of the transversal
variable $y\in\Omega$, and so \eqref{0.1} becomes
\begin{equation} \lb{0.2}
T_t = T_{xx}  + f(T)
\end{equation}
with $x\in\bbR$. This equation has been extensively studied in
mathematical, physical and other literature, starting with the
pioneering works of Fisher \cite{Fisher} and Kolmogorov, Petrovskii,
Piskunov \cite{KPP}. In these papers \eqref{0.2} was used to
describe the propagation of advantageous genes in a population. The
main object of study in these and many subsequent works was the
existence and stability of traveling fronts for \eqref{0.2} and
\eqref{0.1}. In the recent years most of the results have been
extended to include an advection term $u\cdot\nabla T$ in
\eqref{0.1}, and we refer to the reviews \cite{Berrev,Xin2} for an
extensive bibliography.

The above equations are used to model not only population genetics
phenomena. When $f(\tht)>0$ for $\tht\in(0,1)$, then $f$ is a {\it
combustion non-linearity} and \eqref{0.1}/\eqref{0.2} model an
exotermic chemical reaction in an infinite tube with a zero
heat-loss boundary, in particular, flame propagation in a premixed
combustible gas without advection (see Zel'dovich and
Frank-Kamenetskii \cite{ZFK}). In this setting $T$ is the normalized
temperature taking values in $[0,1]$. We note that \eqref{0.1} is
usually obtained from a system involving both the temperature and
the concentration of the reactants after the simplifying assumption
of equal thermal and material diffusivities.

A special case of positive $f$, used often in chemical and
biological literature, is the {\it KPP type} with $f''(\tht)\le c<0$
\cite{KPP}. In combustion models the non-linearity is often
considered to be of {\it Arrhenius type} with slow reaction rates at
low temperatures, modeled by $f(\tht)=e^{-A/\tht}(1-\tht)$. One
often approximates this situation by considering an {\it ignition
non-linearity} $f$ satisfying $f(\tht)=0$ for $\tht\in[0,\tht_0]$
and $f(\tht)>0$ for $\tht\in(\tht_0,1)$, with $\tht_0\in(0,1)$ the
{\it ignition temperature}.

The third prominent case is the {\it bistable non-linearity} with
$f(\tht)< 0$ for $\tht\in(0,\tht_0)$ and $f(\tht)>0$ for
$\tht\in(\tht_0,1)$, where one usually assumes $\int_0^1
f(\tht)\,d\tht>0$. This has been used to model signal propagation
along bistable transmission lines, in particular, nerve pulse
propagation \cite{NAY}. In biological context it is also called {\it
heterozygote inferior} (see Aronson and Weinberger \cite{AW2}).

In this paper we will consider all the above types. Our interest
here will not be in the question of traveling fronts, but in {\it
extinction} of reaction --- {\it quenching} of flames. We will
therefore assume the initial datum $T_0(x)$ for \eqref{0.2} to be
compactly supported, and will want to know when
\begin{equation} \lb{0.3}
\|T(t,\cdot)\|_\infty \to 0 \text{ as } t\to\infty.
\end{equation}
For the sake of simplicity we will restrict ourselves to the case of
$T_0$ being the characteristic function of an interval,
\begin{equation} \lb{0.4}
T_0(x)\equiv \chi_{[-L,L]}(x),
\end{equation}
and study how long-time behavior of $T$ depends on $L$. The methods
in this paper allow one to treat some other increasing 1-parameter
families of initial conditions, too.

Thus, we will study the competition of reaction and diffusion. The
former helps increasing the temperature, whereas the latter
(together with the compactness of the support of the initial datum)
works towards the extinction of the flame. This question was
originally addressed forty years ago by Kanel' \cite{Kanel} who
considered the case of ignition non-linearity and proved that if the
initial datum is large enough, then reaction wins, whereas if it is
small then diffusion manages to quench the flame. More precisely,
when $T$ solves \eqref{0.2}/\eqref{0.4} and $f$ is of ignition type,
Kanel' proved that there are two length scales $L_0$, $L_1$ such
that
\begin{align*}
& T(t,x)\to 0 \text{ as $t\to \infty$ uniformly in $x\in\bbR$ if
$L<L_0$,}
\\ & T(t,x)\to 1 \text{ as $t\to \infty$ uniformly on compacts if
$L>L_1$.}
\end{align*}
This has been extended to the case of bistable $f$ by
Aronson-Weinberger \cite{AW2}. Both results actually hold when
\eqref{0.4} is replaced by
\begin{equation} \lb{0.5}
T_0(x)\equiv \alpha\chi_{[-L,L]}(x).
\end{equation}
for any $\alpha>\tht_0$, with $L_0$ and $L_1$ depending on $\alpha$
(in the ignition case this follows from \cite{Kanel}, in the
bistable case it was proved by Fife and McLeod \cite{FM}). A natural
question arises: does $L_0$ equal $L_1$? If this is true and if one
could determine the behavior of $T$ as $t\to\infty$ when $L=L_0$,
then one would be able to provide the complete ``phase portrait''
for the PDE \eqref{0.2} with respect to a 1-parameter family of
initial conditions.

Since these early works, particularly in the recent years, several
authors have studied quenching for \eqref{0.1}. The above results
have been extended to the case when \eqref{0.1} includes an
advection term $u\cdot\nabla T$, with $u$ a shear or periodic flow
(see \cite{Roq, Xin}), even for certain combustion non-linearities
\cite{ZlaArrh}. Quenching of large initial data by large amplitude
shear and cellular flows has been studied in
\cite{CKR,FKR,KZ,ZlaArrh}. However, the question whether $L_0=L_1$
remained open even in the simplest case of \eqref{0.2}. The
following two results provide the answer, including the treatment of
the critical case $L=L_0$.

The first of them holds for ignition and combustion non-linearities.

\begin{theorem} \lb{T.1.1}
Let $\tht_0\in[0,1)$ and $f:[0,1]\to\bbR$ be Lipschitz with
$f(\tht)=0$ when $\tht\in[0,\tht_0]$, $f(\tht)>0$ when
$\tht\in(\tht_0,1)$, and $f(1)=0$. If $\tht_0>0$ then assume in
addition that $f$ is non-decreasing on $[\tht_0,\tht_0+\del]$ for
some $\del>0$. Let $T:[0,\infty)\times\bbR\to[0,1]$ solve
\begin{align}
& T_t = T_{xx}  + f(T) \lb{1.1}
\\ & T_0(x) \equiv \chi_{[-L,L]}(x). \notag
\end{align}
Then there is $L_0\ge 0$ such that
\begin{SL}
\item[{\rm{(i)}}] if $L<L_0$, then $T\to 0$ uniformly on $\bbR$ as $t\to
\infty$;
\item[{\rm{(ii)}}] if $L=L_0$, then $T\to \tht_0$ uniformly on compacts as $t\to
\infty$;
\item[{\rm{(iii)}}] if $L>L_0$, then $T\to 1$ uniformly on compacts as $t\to
\infty$.
\end{SL}
\end{theorem}

{\it Remark.} The possibility of $L_0=0$ (so called {\it
hair-trigger effect}) cannot be excluded when $\tht_0=0$. More
precisely, using results from \cite{ZlaArrh} one can show that
$L_0=0$ when $f(\tht)\ge c\tht^p$ for some $p<3$ and all small
$\tht$, but $L_0>0$ when $f(\tht)\le c\tht^p$ for some $p>3$ and all
small $\tht$. We also note that if $\tht_0=0$, then the convergence
in (ii) is as in (i)
--- uniform on $\bbR$.
\smallskip

Our second result holds for bistable non-linearities. We define
$\tht_2\in(\tht_0,1)$ by $\int_0^{\tht_2} f(\tht)d\tht=0$ and let
$U$ be the unique function solving $0=U''+f(U)$ with $U(0)=\tht_2$
and $U'(0)=0$. Then $U$ is an even function and we will show in the
proof of the following theorem that it is positive on $\bbR$,
decreasing to 0 on $(0,\infty)$, and bell-shaped.

\begin{theorem} \lb{T.1.2}
Let $\tht_0\in(0,1)$ and $f:[0,1]\to\bbR$ be Lipschitz with
$f(0)=f(\tht_0)=f(1)=0$, $f(\tht)< 0$ when $\tht\in(0,\tht_0)$, and
$f(\tht)>0$ when $\tht\in(\tht_0,1)$. Assume also that $\int_0^1
f(\tht)d\tht>0$. Let $T:[0,\infty)\times\bbR\to[0,1]$ solve the
problem \eqref{1.1}.
Then there is $L_0>0$ such that
\begin{SL}
\item[{\rm{(i)}}] if $L<L_0$, then $T\to 0$ uniformly on $\bbR$ as $t\to
\infty$;
\item[{\rm{(ii)}}] if $L=L_0$, then $T\to U$ uniformly on $\bbR$ as $t\to
\infty$;
\item[{\rm{(iii)}}] if $L>L_0$, then $T\to 1$ uniformly on compacts as $t\to
\infty$.
\end{SL}
\end{theorem}

{\it Remark.} Both results can be extended to some other increasing
families of initial conditions. In particular, to \eqref{0.5} with
$\alpha>\tht_0$.
\smallskip


The crux of the proofs of both theorems will be to show that there
is a single $L$ for which $T$ does not converge to 0 or 1 at $x=0$
as $t\to\infty$. In Theorem \ref{T.1.1} this will be achieved with
the help of Lemma~\ref{L.2.2} by comparing solutions of \eqref{1.1}
for two different initial conditions at differently rescaled times.
In Theorem \ref{T.1.2} it will follow from a detailed analysis of
the large time behavior of $T$ when the above limit is not 0 or 1,
and an application of the comparison principle.

We note here that Theorem \ref{T.1.1} is, in a sense, a limiting
case of Theorem \ref{T.1.2}. If one takes $f\to 0$ on $(0,\tht_0)$
keeping $f$ unchanged on $(\tht_0,1)$, one has $\tht_2\to\tht_0$ and
$U\to\tht_0$ on compacts. That is, the bell shaped solution $U$ from
Theorem \ref{T.1.2}(ii) converges to the constant solution $\tht_0$
from Theorem \ref{T.1.1}(ii).

The rest of the paper is devoted to the proofs of the two theorems.
Section \ref{S2} contains preliminary Lemmas \ref{L.2.1} and
\ref{L.2.2}. Sections \ref{S3} and \ref{S4} prove Theorems
\ref{T.1.1} and \ref{T.1.2}, respectively.

The author thanks Henri Berestycki, Peter Constantin, Fran\c cois
Hamel, Alexander Kiselev, Peter Pol\' a\v cik, Jean-Michel
Roquejoffre, and Lenya Ryzhik for encouragement and useful
discussions. He also acknowledges the hospitality of the Mathematics
Department of the University of Chicago where part of this work was
done, as well as partial support from the NSF through the grant
DMS-0314129.

\section{Preliminary Lemmas} \lb{S2}

We start with

\begin{lemma} \lb{L.2.1}
Let $f:[0,1]\to\bbR$ be Lipschitz with $f(0)=f(1)=0$. If
$T:[0,\infty)\times\bbR\to[0,1]$ solves \eqref{1.1},
then the following hold.
\begin{SL}
\item[{\rm{(i)}}] If $|x|\le|y|$ then $T(t,x)\ge T(t,y)$.
\item[{\rm{(ii)}}] There is $t_\ast>0$ (possibly $t_\ast=\infty$)
such that $T(t,0)$ as a function of $t$ is non-increasing on
$[0,t_\ast)$ and non-decreasing on $[t_\ast,\infty)$.
\item[{\rm{(iii)}}] If $f\ge 0$ then there is $\tht_*\in[0,1]$ such
that $f(\tht_*)=0$ and $T(t,x)\to\tht_*$ as $t\to\infty$, uniformly
on compacts.
\end{SL}
\end{lemma}

{\it Remarks.} 1. For sufficiently smooth $f$ this is essentially a
result of Kanel' \cite{Kanel}, although in (ii) he only proves
$T(t,0)$ to be eventually monotone.
\smallskip

2. In the case of \eqref{0.5} with $\alpha\in(\tht_0,1)$, part (ii)
has $0<t_\ast\le t_{\ast\ast}\le\infty$ such that $T(t,0)$ is
non-decreasing on $[0,t_\ast)$, non-increasing on
$[t_\ast,t_{\ast\ast})$ and non-decreasing on
$[t_{\ast\ast},\infty)$.
\smallskip

\begin{proof}
We first assume that $f$ is smooth and briefly recall main points of
the proofs of (i) and (ii) from \cite{Kanel}. Let $T^\eps$ solve
\eqref{1.1} but with initial condition $T(0,x)\equiv \chi^\eps (x)$
where $\chi^\eps$ are smooth, symmetric, decreasing in $|x|$, and
converge to $\chi_{[-L,L]}$ in $L^1(\bbR)$ as $\eps\to 0$. Then
$T^\eps_x(0,x)\le 0$ for $x>0$, and by symmetry $T^\eps_x(t,0)=0$.
Since
\[
(T^\eps_x)_t = (T^\eps_x)_{xx}  + f'(T^\eps)T^\eps_x,
\]
the maximum principle gives $T^\eps_x(t,x)\le 0$ for $x>0$. Symmetry
then yields $T^\eps_x(t,x)\ge 0$ for $x<0$. Since for any fixed
$t>0$ we have $T^\eps(t,x)\to T(t,x)$ uniformly in $x$ as $\eps\to
0$, this proves (i).

Now let $D^h(t,x)\equiv T(t+h,x)-T(t,x)$. By the mean value theorem,
\[
D^h_t = D^h_{xx} + f'(S)D^h
\]
for some $S(t,x)$. Let $\Delta_h$ be set of $(t,x)$ for which
$D_h(t,x)\le 0$. Then $\Delta_h\cap (\{0\}\times\bbR) =
\{0\}\times[-L,L]$. By the maximum principle and symmetry,
$\Delta_h$ is connected and its sections by lines parallel to the
$x$-axis are segments symmetric about the $t$-axis. Therefore there
is $0<t^h_*\le\infty$ such that $D^h(t,0)<0$ for $t\in[0,t^h_*)$ and
$D^h(t,0)\ge 0$ for $t\in[t^h_*,\infty)$. From $D^h(t,x) =
D^{h/2}(t+\tfrac h 2) + D^{h/2}(t)$ we obtain $t^{h/2}_*
\in[t^h_*,t^h_*+\tfrac h2]$, and (ii) follows with $t_*\equiv
\lim_{n\to\infty} t^{2^{-n}}_*$.

If $f$ is only Lipschitz, take smooth $f^\eps$ such that $\|
f^\eps-f\|_\infty \le \eps$ and let $T^\eps$ solve \eqref{1.1} with
$f^\eps$ in place of $f$. One can then show that $V^\eps\equiv
T^\eps -T$ satisfies $|V^\eps(t,x)|\le \tfrac \eps c(e^{ct}-1)$ with
$c$ the Lipschitz constant for $f$ (we spell this argument out in
the proof of Theorem \ref{T.1.1} below). Therefore
$T(t,x)=\lim_{\eps\to 0} T^\eps(t,x)$ for all $t$ and $x$, and since
(i) and (ii) hold for each $T^\eps$, they also hold for $T$.

Finally, assume that $f\ge 0$. By (ii), $\tht_*\equiv
\lim_{t\to\infty} T(t,0)$ is well defined. Let $\Phi$ solve
$\Phi_t=\Phi_{xx}$ on $\bbR^+$ with $\Phi(0,x)\equiv T(0,x)$ and
boundary condition $\Phi(t,0)\equiv T(t,0)$. Then
$\Phi(t,x)\to\tht_*$ as $t\to\infty$, uniformly on compacts. Since
by the comparison principle (e.g., \cite{Sm}) and (i), $\Phi(t,x)\le
T(t,x)\le T(t,0)$, the second claim in (iii) follows.

To prove the first claim, assume $f(\tht_*)>0$ and choose $\eps>0$
such that for $\tht\le\tht_*+10\eps$ we have $f(\tht)\ge
\tht-\tht_*+2\eps$. Pick $t_0$ such that if $\Phi$ solves
$\Phi_t=\Phi_{xx}$ on $\bbR$ with initial condition
$\Phi(t_0,x)=T(t_0,x)$, then $\Phi(t,0)\ge \tht_*-\eps$ and
$T(t,x)\le \tht_*+\eps$ for $t\in[t_0,t_0+\ln 4]$ and $x\in\bbR$.
This is possible thanks to the second claim in (iii). Define
\[
S(t,x)\equiv \tht_*-2\eps + (\Phi(t,x)-\tht_*+2\eps)e^{t-t_0}.
\]
Then $S(t,x)\le \tht_*+10\eps$ for $t\in[t_0,t_0+\ln 4]$ because
$\Phi(t,x)\le T(t,x)\le\tht_*+\eps$ for these $t$. A simple
computation now shows that $S_t\le S_{xx} + f(S)$ for
$t\in[t_0,t_0+\ln 4]$. Hence, $S$ is a subsolution of \eqref{1.1}
with $S(t_0,x)=T(t_0,x)$, and so $S\le T$ for $t\in[t_0,t_0+\ln 4]$.
But $S(t_0+\ln 4,0)\ge \tht_*+2\eps>T(t,0)$, which is a
contradiction. Therefore we have to have $f(\tht_*)=0$.
\end{proof}

Next, observe that we can use scaling to replace the variation in
the initial condition in \eqref{1.1} by variation in the reaction
strength. If $T$ solves \eqref{1.1} with $T(0,x)\equiv
\chi_{[-L,L]}(x)$, define $\til T(t,x) \equiv T(L^2 t,Lx)$, so that
we have
\[
\til T_t = \til T_{xx} + L^2f(\til T)
\]
and $\til T(0,x)=\chi_{[-1,1]}(x)$. Hence, Theorem \ref{T.1.1} will
be proved if we show that its conclusion holds for the $L$-dependent
family of problems
\begin{align}
& T_t  = T_{xx}  + Lf(T) \lb{2.11}
\\ & T_0(x)  \equiv \chi_{[-1,1]}(x) \notag
\end{align}
instead of \eqref{1.1} (note that Lemma \ref{L.2.1} holds here,
too). This important observation motivates the following key lemma.

\begin{lemma} \lb{L.2.2}
Let $\Omega$ be a connected open domain in $\bbR^n$ with a smooth
boundary (possibly $\Omega=\bbR^n$)
and let $f,g:[0,\infty)\to\bbR$ be Lipschitz with $f(0)=g(0)=0$ and
$f\le g$. Let $T,S:[0,\infty)\times\Omega\to[0,\infty)$ be
continuous functions solving
\begin{align}
T_t &=\Delta T  + f(T) \lb{2.1}
\\ S_t &=\Delta S  + g(S) \lb{2.2}
\end{align}
in $\Omega$ with Dirichlet boundary conditions on $\partial\Omega$.
Assume $0\le T(0,x)\le S(0,x)$ for all $x\in\Omega$ and
$T(0,x_0)<S(0,x_0)$ for some $x_0$. Assume also that for any
$\tht>0$ the set $\Omega_{0,\tht}\equiv \{ x\in\Omega \,|\,
S(0,x)\ge\tht \}$ is compact. Finally, assume that there are
$\tht_1>0$ and $\eps_1>0$ such that for any
$\tht\in[\tht_1,\|T\|_\infty)$ and $\eps\in[0,\eps_1]$ we have
\begin{equation} \lb{2.3}
g \big( \tht + \eps [\tht-\tht_1] \big) \ge (1+\eps) f (\tht).
\end{equation}
Then
\begin{equation} \lb{2.4}
\liminf_{t\to\infty} \inf_{T(t,x)>\tht_1} \frac {S(t,x)-\tht_1}
{T(t,x)-\tht_1} > 1
\end{equation}
with the convention that infimum over an empty set is $\infty$.
\end{lemma}

{\it Remark.} The result holds without change when we add a first
order term $u(x)\cdot\nabla$ with $u\in C^1$ to \eqref{2.1} and
\eqref{2.2}.

\begin{proof}
First notice that the assumptions imply that
\begin{equation*}
\Omega_{t,\tht}\equiv \{ x\in\Omega \,\big|\, S(t,x)\ge\tht \}
\end{equation*}
are compact. Indeed, by the maximum principle,
$\Omega_{t,\tht}\subseteq\til\Omega_{t,\del\tht}$ where $\del\equiv
e^{-tc}$ with $c$ the Lipschitz constant for $g$, and
$\til\Omega_{t,\tht}$ is defined as $\Omega_{t,\tht}$ but with
$\Phi$, the solution of
\[
\Phi_t = \Delta\Phi, \quad\quad \Phi(0,x)=S(0,x),
\]
in place of $S$. Compactness of $\til\Omega_{t,\tht}$ follows from
that of $\til\Omega_{0,\tht}$ and the Feynman-Kac formula.

The assumptions and the strong maximum principle also imply
$T(t,x)<S(t,x)$ for $t>0$ and $x\in\Omega$.
Let us define
\begin{align*}
\Omega(t) &\equiv \{ x\in\Omega \,\big|\, T(t,x)>\tht_1 \},
\\ \Omega'(t) &\equiv \{ x\in\Omega \,\big|\, T(t,x)=\tht_1 \},
\end{align*}
and let
\[
\omega(t)\equiv \min \bigg\{ 1+\eps_1, \inf_{x\in\Omega(t)} \frac
{S(t,x)-\tht_1} {T(t,x)-\tht_1} \bigg\}.
\]
Since $\overline{\Omega(t)}$ is compact and $T<S$ continuous,
$\omega(t)>1$ for $t>0$. Hence the result will follow if we show
that $\omega$ is a non-increasing function. Since $\omega$ is
continuous (because $T$ and $S$ are), it is sufficient to show that
for any $t_0>0$ there is $\tau_0>0$ such that for all
$t\in[t_0,t_0+\tau_0]$ we have $\omega(t)\ge\omega(t_0)$.

Hence, fix $t_0>0$. Notice that uniform boundedness of the sets
$\Omega_{t,\tht_1}$ for $t\in[0,t_0+1]$ (by the maximum principle
and Feynman-Kac formula) and continuity of $T,S$ imply that $T,S$
are uniformly continuous on $[0,t_0+1]\times\Omega$, and the set
\[
\Sigma\equiv \{ (t,x) \,|\, t\in[t_0,t_0+1] \text{ and }
x\in\Omega(t)\cup\Omega'(t) \}
\]
is compact. Thanks to the uniform continuity of $T$ we only need to
consider the case $\Omega(t_0)\cup\Omega'(t_0)\neq\emptyset$, and
hence $\Sigma\neq\emptyset$. Since $S$ is continuous and $S>T$,
\[
\sigma\equiv \inf_{(t,x)\in\Sigma} \{S(t,x)-\tht_1 \}>0.
\]
We let $\del\equiv \sigma/4(1+\eps_1)$ and define
\begin{equation*}
\Delta \equiv \{ x \,\big|\, |T(t_0,x)-\tht_1|\le \del \text{ and }
S(t_0,x)-\tht_1\ge \sigma-\del \}
\end{equation*}
By the uniform continuity of $T,S$, there is $\tau_0\in(0,1)$ such
that for $t\in[t_0,t_0+\tau_0]$ and $x\in\Omega$ we have
\begin{equation} \lb{2.5}
|T(t,x)-T(t_0,x)|\le\tfrac\del 2 \quad \text{ and } \quad
|S(t,x)-S(t_0,x)|\le\tfrac\del 2.
\end{equation}
So if $t\in[t_0,t_0+\tau_0]$, then
\begin{equation} \lb{2.7}
\Omega(t)\subseteq \Omega(t_0)\cup\Delta
\end{equation}
(note that $S(t_0,x)-\tht_1\ge \sigma-\tfrac\del 2$ for
$x\in\Omega(t)$ because then $(t,x)\in\Sigma$). Now if
$t\in[t_0,t_0+\tau_0]$ and $x\in\Delta$, then by \eqref{2.5},
\begin{equation} \lb{2.6}
S(t,x)-\tht_1 > \frac \sigma 2 > (1+\eps_1)|T(t,x)-\tht_1| \ge
\omega(t_0) (T(t,x)-\tht_1).
\end{equation}

Next let
\begin{equation} \lb{2.8}
A\equiv  \{ x\in\Omega \,\big|\, T(t_0,x)>\tht_1+\del \} =
\Omega(t_0)\smallsetminus\Delta
\end{equation}
and
\[
B\equiv \{ x\in\Omega \,\big|\, T(t_0,x)\ge\tht_1+\tfrac\del 2 \}
\subseteq\Omega(t_0).
\]
Uniform continuity of $T$ shows that, $\dist(A,B^c)>0$, and so there
is an open set $\Gamma$ with a smooth boundary such that
$A\subseteq\Gamma\subseteq B$. Let $\til T\equiv T-\tht_1$, $\til
U\equiv \omega(t_0) \til T$, $\til S\equiv S-\tht_1$, $\til
f(\tht)\equiv f(\tht+\tht_1)$, and $\til g(\tht)\equiv
g(\tht+\tht_1)$. Then for $x\in\Gamma$ we have
\[
\til S(t_0,x)\ge \omega(t_0) \til T(t_0,x)=\til U(t_0,x)
\]
by the definition of $\omega(t_0)$, for $t\in[t_0,t_0+\tau_0]$ and
$x\in\partial\Gamma$ we have
\[
\til S(t,x) > \sigma-2\del > \frac \sigma 2 \ge \omega(t_0) \til
T(t,x)=\til U(t,x)
\]
since $\partial\Gamma\subseteq B\smallsetminus A\subseteq\Delta$,
and for $t\in[t_0,t_0+\tau_0]$ and $x\in\Gamma$ we have
\begin{align*}
\til U_t &=\Delta\til U + \omega(t_0) \til f \big( \tfrac
1{\omega(t_0)} \til U \big),
\\ \til S_t &=\Delta\til S + \til g(\til S)
\end{align*}
by \eqref{2.1} and \eqref{2.2}. For these $(t,x)$ we have
$T(t,x)\ge\tht_1$ because of \eqref{2.5} and $\Gamma\subseteq B$,
and so by \eqref{2.3} and $\omega(t_0)-1 \in(0,\eps_1]$,
\[
\omega(t_0) \til f \big( \tfrac 1{\omega(t_0)} \til U \big) =
\omega(t_0) f(T)\le g\big( \omega(t_0)[T-\tht_1]+\tht_1 \big)=\til
g(\til U).
\]
The comparison principle now shows that $\til S\ge\til U$ on
$[t_0,t_0+\tau_0]\times\Gamma$. Hence
\[
S(t,x)-\tht_1 \ge \omega(t_0) (T(t,x)-\tht_1)
\]
for $t\in[t_0,t_0+\tau_0]$ and $x\in A$, which together with
\eqref{2.6}, \eqref{2.8}, and \eqref{2.7} gives
$\omega(t)\ge\omega(t_0)$ for $t\in [t_0,t_0+\tau_0]$. The proof is
finished.
\end{proof}

\section{Proof of Theorem \ref{T.1.1}} \lb{S3}

We can now complete the proof of Theorem \ref{T.1.1}. We will do
this for the formulation in \eqref{2.11}.

First assume $\tht_0>0$. We know from Lemma \ref{L.2.1}(iii) that
for every $L$ we have $T\to\tht^L_*$ uniformly on compacts, with
$\tht_*^L$ such that $f(\tht^L_*)=0$. Obviously
$\tht^L_*\notin(0,\tht_0)$ because in that case we would have
$T(t,x)\le\tht_0$ for all $t\ge t_0$ and consequently $T\to 0$
(since $\|T(t_0,\cdot)\|_1 <\infty$ and $T_t=T_{xx}$ for $t\ge
t_0$). So we are only left with $\tht^L_*\in\{0,\tht_0,1\}$.

Let $A$, $B$, and $C$ be the sets of $L\ge 0$ such that $\tht^L_*$
equals $0$, $\tht_0$, and $1$, respectively. Notice that since
$T(t,0)\ge T(t,x)$, the convergence of $T$ to $0$ for $L\in A$ is
actually uniform on $\bbR$. We have $A\cup B\cup C=[0,\infty)$ and
the comparison principle implies that the three sets are intervals
with $A$ lying to the left of $B$ and $B$ to the left of $C$.

Moreover, $A$ and $C$ are non-empty by Kanel' \cite{Kanel} and open.
The latter follows from the fact that if $T^L$ is the solution of
\eqref{2.11}, then for $L_1<L_2$ and $V\equiv T^{L_2}-T^{L_1}$ we
have $V\ge 0$ by comparison, and
\begin{align*}
V_t & = \Delta V +(L_2-L_1)f(T^{L_2}) + L_1[ f(T^{L_2})-f(T^{L_1}) ]
\\ & \le \Delta V + c(L_2-L_1) + cL_1 V
\end{align*}
with $c\ge\|f\|_\infty$ the Lipschitz constant for $f$. Since the
function $\til V(t,x)\equiv \tfrac{L_2-L_1}{L_1}(e^{cL_1 t}-1)$
satisfies
\[
\til V = \Delta \til V + c(L_2-L_1) + cL_1 \til V
\]
with $\til V(0)=0=V(0)$, the comparison principle gives $V\le\til
V$, that is,
\[
T^{L_2}(t,x)-T^{L_1}(t,x) \in [0, \tfrac{L_2-L_1}{L_1}(e^{cL_1
t}-1)].
\]
Therefore if $L_1\in A$, then $T^{L_1}(t_0,0)\le\tfrac 12 \tht_0$
for some $t_0>0$, and hence $T^{L_2}(t_0,0)<\tht_0$ (and so $L_2\in
A$) for $L_2<L_1+ \tfrac 12 \tht_0 L_1 (e^{cL_1 t_0}-1)^{-1}$. On
the other hand, Kanel's result \cite{Kanel} also holds for
\eqref{0.5}, and it says that for any $\alpha>\tht_0$ and $L>0$
there is $M=M(\alpha,L)<\infty$ such that if $T$ solves \eqref{2.11}
and $T(t_0,x)\ge \alpha\chi_{[-M,M]}(x)$, then $T\to 1$ uniformly on
compacts. Let $\tht_0<\alpha<\beta<1$ and if $L_1\in C$, let
$M=M(\alpha,\tfrac 12 L_1)$. For some $t_0$ we have
$T^{L_1}(t_0,x)\ge\beta\chi_{[-M,M]}(x)$ and so for any
$L_2>L_1-(\beta-\alpha)L_1(e^{cL_1 t_0}-1)^{-1}$ we have
$T^{L_2}(t_0,x)\ge\alpha\chi_{[-M,M]}(x)$. If in addition
$L_2>\tfrac 12 L_1$, we have $L_2\in C$. So $A,C$ are non-empty and
open, and hence $B$ is non-empty and closed.

The proof will be finished if we show that $B$ contains a single
element. Hence assume $L_1<L_2$ are both in $B$. Let $\tht_1\equiv
\tfrac 12 \tht_0 \in(0,\tht_0)$ and
\[
\eps_1\equiv \min\{ L_2 L_1^{-1} -1,\, \del(\del+\tht_0)^{-1} \} >0
\]
(with $\del$ from the statement of Theorem \ref{T.1.1}). Choose
$t_0>0$ such that $T^{L_1}(t,x)\le \tht_0+\tfrac\del 2$ when $t\ge
t_0$. The comparison principle, $f\not\equiv 0$, and the strong
maximum principle yield $T^{L_1}< T^{L_2}$ for $t>0$ and both
$T^{L_1}$ and $T^{L_2}$ are obviously continuous for $t>0$.
Lipschitzness of $f$ and compact support of $T^{L_2}(0,\cdot)$ show
that for any $\tht>0$, the set of $x$ for which
$T^{L_2}(t_0,x)\ge\tht$ is compact. Finally, whenever
$\tht\in[\tht_1,\tht_0+\tfrac\del 2]$ and $\eps\in[0,\eps_1]$, we
have $\tht+\eps[\tht-\tht_1]\le \tht_0+\del$. Thus by the
assumptions on $f$ (and the definition of $\eps_1$),
\[
L_2 f(\tht+\eps[\tht-\tht_1])\ge L_2 f(\tht) \ge (1+\eps)L_1
f(\tht).
\]
Therefore Lemma \ref{L.2.2} applies to $T^{L_1}$ and $T^{L_2}$ (with
starting time $t_0$) and shows that for some $r>1$ and all large
enough $t$ we have
\[
T^{L_2}(t,x)-\tht_1 \ge r\big[ T^{L_1}(t,x)-\tht_1 \big]
\]
whenever $T^{L_1}(t,x)>\tht_1$. But since $\tht_1<\tht_0$, this
contradicts the assumption that both $T^{L_1}(t,0)$ and
$T^{L_2}(t,0)$ converge to $\tht_0$ as $t\to\infty$. Hence, $B=\{
L_0\}$ and we are done.

Now, consider $\tht_0=0$. We have $\tht^L_*\in\{0,1\}$, the sets
$A,C$ satisfy $A\cup C=[0,\infty)$, and by the comparison principle,
$A$ lies to the left of $C$. Moreover, $0\in A$ and $C$ is non-empty
and open by the same argument as above. Hence $A$ is closed and its
maximum is $L_0$ (possibly $L_0=0$). Lemma \ref{L.2.1}(iii) yields
(iii) of this theorem and $T(t,0)\ge T(t,x)$ gives (i) and (ii),
including the fact that the convergence in (ii) is uniform on
$\bbR$. The proof is finished.

\section{Proof of Theorem \ref{T.1.2}} \lb{S4}

The situation is somewhat more complicated here. Firstly, we do not
have Lemma \ref{L.2.1}(iii) at our disposal, and so the limit of $T$
as $t\to\infty$ need not always be a constant function. And
secondly, we cannot use Lemma \ref{L.2.2} and the scaling argument
preceding it in the way we did in the last section because it is not
anymore true that $L_2f\ge L_1f$ when $L_2>L_1$. We note that one
can still use the lemma without scaling, but then the argument
applies only to a restricted class of bistable $f$. Fortunately, it
turns out that the first of these difficulties actually cancels the
problems created by the second, as we shall see below.


Let us therefore go back to $T$ solving \eqref{1.1} rather than
\eqref{2.11}. We know from Lemma \ref{L.2.1}(ii) that
$\tht^L_*\equiv\lim_{t\to\infty} T(t,0)$ is well defined, and from
the comparison principle that it is non-decreasing in $L$.

First assume $\tht^L_*<\tht_2$, with $\tht_2$ defined in the
introduction. Choose $\eps>0$ and a Lipschitz function $\til
f:[0,1]\to\bbR$ so that $\til f=0$ on $[0,\eps]$, $\til f'(\eps)<0$,
$\til f\ge f$ on $(\eps,\tfrac 12(\tht^L_*+\tht_2)]$ and $\til f$
has a single zero there, $\til f>0$ on $(\tfrac
12(\tht^L_*+\tht_2),1)$, $\til f(1)=0>\til f'(1)$, and
\begin{equation} \lb{3.1}
\int_0^1 \til f(\tht) d\tht<0.
\end{equation}
Let $t_0$ be such that for $t\ge t_0$ and all $x\in\bbR$ we have
$T(t,x)\le \tfrac 12(\tht^L_*+\tht_2)$. This is possible by Lemma
\ref{L.2.1}(i). Since $\til f\ge f$ on $[0,\tfrac
12(\tht^L_*+\tht_2)]$, starting from time $t_0$ one has $T_t\le
T_{xx}+\til f(T)$, that is, $T$ is a subsolution of the equation
\begin{equation} \lb{3.2}
\Phi_t= \Phi_{xx}+\til f(\Phi).
\end{equation}

Let $\phi:\bbR\to[0,1]$ with $\phi(x)\to \eps$ as $x\to\infty$ and
$\phi(x)\to 1$ as $x\to -\infty$ be the unique, up to translation,
traveling front profile (with speed $v$) for \eqref{3.2}
\cite{Kanel2}. That is, $\phi(x-vt)$ solves \eqref{3.2}. It follows
from \eqref{3.1} that in this case $v<0$.

From compactness of the support of $T(0,x)$ and Lipschitzness of
$f$, $T(t_0,x)\to 0$ as $|x|\to\infty$. This and
$\|T(t_0,\cdot)\|_\infty<1$ mean that there is $x_0$ such that
$T(t_0,x)\le\phi(x-x_0-vt_0)$, and since $\phi(x-x_0-vt)$ is a
solution and $T(t,x)$ a subsolution of \eqref{3.2},
\[
T(t,x)\le\phi(x-x_0-vt)
\]
for all $t\ge t_0$. But then $T(t,0)\le \phi(-x_0-vt) \to \eps$ as
$t\to\infty$ because $v<0$. This holds for any small enough $\eps>0$
and thus $\tht^L_*=0$.

Next assume $\tht^L_*>\tht_2$. Let $S$ be the solution of
\eqref{1.1} on $\bbR^+$ with $S(0,x)=0$ and $S(t,0)=s(t)$ a smooth
strictly increasing function with all derivatives bounded such that
$s(0)=0$, $s(t)\le T(t,0)$, and $\lim_{t\to\infty}s(t)=\tht^L_*$.
Then for any $h>0$ we have $S(h,x)> S(0,x)$ and so by comparison
$S(t+h,x)> S(t,x)$. Hence $\til S(x)\equiv \lim_{t\to\infty}
S(t,x)>0$ is well defined and $\til S(0)=\tht^L_*$. Since by
comparison again $S(t,x)\le T(t,x)\le T(t,0)$, we also have $\til
S(x)\le\tht^L_*$.


Standard parabolic regularity shows that $S(t,x)$ converges to $\til
S(x)$ along with its first two derivatives uniformly on compacts,
and so $\til S$ solves the stationary problem
\begin{equation} \lb{3.3}
0=\til S''+f(\til S)
\end{equation}
on $\bbR^+$ (this can be found also in \cite{AW2}). But then for any
$y>0$
\[
\int_{\til S(y)}^{\tht^L_*} f(\tht) \,d\tht = \int_y^0 f(\til
S(x))\til S'(x) \,dx = \int_0^y \til S''(x)\til S'(x) \,dx = \tfrac
12 \big[ (\til S'(y))^2 - (\til S'(0))^2 \big].
\]
Assume there is $y>0$ such that $\til S(y)<\tht^L_*$, and then pick
one such that also $\til S'(y)<0$. Since $\int_{\omega}^{\tht^L_*}
f(\tht) d\tht$ is bounded below by a positive constant for all
$\omega\in[0,\til S(y)]$, so is $-\til S'(z)$ for all $z\ge y$
(because $\til S(z)\ge 0$). But that means $\til S$ is eventually
negative, a contradiction. Hence, we must have $\til
S\equiv\tht^L_*$, which is only possible if $\tht_*^L=1$. Moreover,
since $S$ converges to $\til S\equiv 1$ uniformly on compacts as
$t\to\infty$ (and $S\le T\le 1$), so does $T$.

The above shows that $\tht^L_*\in\{ 0,\tht_2,1\}$. As in the proof
of Theorem~\ref{T.1.1}, and using the equivalent of Kanel's result
for \eqref{0.5} and bistable $f$ \cite{FM}, one can show that the
intervals $A$, $C$ of $L$ for which $\tht_*^L=0,1$, respectively,
are non-empty and open. If $B$ is the interval of $L$ for which
$\tht_*^L=\tht_2$, then $B$ lies between $A$ and $C$ and again
$A\cup B\cup C=[0,\infty)$.

Next we need to prove that $B$ only contains one element.  We will
show below that if $L\in B$, then $T(t,x)\to U(x)$ uniformly on
$\bbR$ as $t\to\infty$. Here $U$ solves \eqref{3.3} with
$U(0)=\tht_2$ and $U'(0)=0$. Assume now that $L_1<L_2$ are both in
$B$, with $T^{L_1}$ and $T^{L_2}$ the corresponding solutions of
\eqref{1.1}. We then have $T^{L_1}(t,0)\to\tht_2$, and since the
equation is translation invariant, we also have $\til
T(t,\eps)\to\tht_2$ when $\til T$ solves \eqref{1.1} with initial
condition $\til T_0(x)\equiv \chi_{[-L_1+\eps,L_1+\eps]}(x)$. But if
$|\eps|<L_2-L_1$, then $T^{L_2}_0(x)\ge \til T_0(x)$, and so by the
comparison principle,
\[
U(\eps) = \lim_{t\to\infty} T^{L_2}(t,\eps)\ge \tht_2.
\]
Since $U''(0)=-f(U(0))=-f(\tht_2)<0$, $U$ has a strict local maximum
at zero and therefore $U(\eps)<U(0)=\tht_2$ for all small enough
$|\eps|>0$. This is a contradiction and hence $B=\{ L_0\}$.

To complete the proof, we need to show that $T(t,x)\to U(x)$
uniformly as $t\to\infty$ when $L\in B$ (and hence
$\tht_*^L=\tht_2$).  Notice that the argument following \eqref{3.3}
applies to $U$ and we find for any $x>0$ such that $U(x)\ge 0$,
\begin{equation} \lb{3.4}
\int^{\tht_2}_{U(x)} f(\tht)\,d\tht = \tfrac 12 (U'(x))^2.
\end{equation}
The definition of $\tht_2$ then shows that $U(x)\le\tht_2$, and
$U'(x)\neq 0$ when $U(x)\in(0,\tht_2)$. Since $U(x)$ cannot be
constant $\tht_2$ on any interval and $U'$ is continuous, we must
have $U'(x)<0$ for all $x>0$ such that $U(x)>0$. There is no $x$
with $U(x)=0$ because then \eqref{3.4} would give $U'(x)=0$,
contradicting uniqueness of solutions to initial value problems
associated to \eqref{3.3}. Hence $U(x)\in(0,\tht_2)$ and $U'(x)<0$
for $x>0$, with $U'(x)$ bounded away from zero when $U(x)$ is away
from zero (by \eqref{3.4} and the definition of $\tht_2$). This and
symmetry show that $U$ is indeed a symmetric bell-shaped solution
(with $U'$ decreasing on $[0,U^{-1}(\tht_0)]$ and increasing on
$[U^{-1}(\tht_0),\infty)$ by \eqref{3.4}) of the stationary problem
\eqref{3.3} such that $U(x)\to 0$ as $|x|\to\infty$.

If we now apply the argument involving $S$ and $\til S$ from the
case $\tht_*^L>\tht_2$, we find as above that $\til
S(x)\le\tht_2=\til S(0)$ and that $\til S>0$ is possible only if
$\til S'(0)= 0$. But then $\til S(0)=U(0)$ and $\til S'(0)=U'(0)$,
thus $\til S=U$. Moreover, uniform on compacts convergence of $S$ to
$U$ and $0\le S(t,x)\le U(x)\to 0$ as $|x|\to\infty$ yield uniform
on $\bbR$ convergence of $S$ to $U$. Since $T(t,x)\ge S(t,x)$, we
have $\liminf_{t\to\infty} T(t,x)\ge U(x)$ uniformly on $\bbR$.
Hence we are left with proving $\limsup_{t\to\infty} T(t,x)\le U(x)$
uniformly in $x>0$ (which suffices due to symmetry).

Let $0<x_0<\infty$ be such that if $S(0,x)\ge \tht_2
\chi_{[-x_0,x_0]}(x)$ and $S$ satisfies \eqref{1.1}, then $S\to 1$
uniformly on compacts. Such $x_0$ exists by \cite{FM} because
$\tht_2>\tht_0$. Then obviously for every $t\ge 0$ we have
$T(t,x_0)\le\tht_2$, because otherwise Lemma \ref{L.2.1}(i) would
imply $T\to 1$ uniformly on compacts. Since both $T(t,x)$ and
$V(x)=U(x-x_0)$ satisfy \eqref{1.1} on $(x_0,\infty)$,
$V(x_0)=\tht_2\ge T(t,x_0)$, and $V(x)>0=T(0,x)$ for $x>x_0$, the
comparison principle implies $T(t,x)\le V(x)$. Let us therefore
define
\[
x_1\equiv \min\{\til x\,|\, \limsup_{t\to\infty} T(t,x)\le U(x-\til
x) \text{ uniformly in $x>\til x$}\} \le x_0.
\]
The minimum is achieved because $U$ is uniformly continuous and by
Lemma \ref{L.2.1}(i), $\limsup_{t\to\infty} T(t,x)\le\tht_2=U(0)$.
We note that at this point one can derive $T\to U$ from $x_1<\infty$
and the results of \cite{BJP} when $f\in C^1$ and $f'(0)<0$. Our
non-linearity is more general and so \cite{BJP} is not applicable
here.

If $x_1=0$, then we are done, so assume $x_1>0$. First notice that
$\limsup_{t\to\infty} T(t,\tfrac 12 x_1)\le \tht_2-\del_1$ for some
$\del_1>0$. Indeed
--- in view of $\lim_{t\to\infty} T(t,0)=\tht_2$, Lemma
\ref{L.2.1}(i), and the comparison principle, it is sufficient to
show that there are $\del_1,\del_2>0$ such that if $S(0,x)\ge
(\tht_2-\del_1)\chi_{[-x_1/2,x_1/2]}(x)$ and $S$ satisfies
\eqref{1.1}, then $S(t,0)\ge \tht_2+\del_2$ for some $t>0$. This in
turn is true because it holds for $\del_1=0$ and some $t,\del_2>0$,
since then $S(0,0)=\tht_2$ and $S_t(0,0)=f(\tht_2)>0$, and because
$S(t,0)$ is continuous in $\del_1$.

Now choose $x_2\in(\tfrac 12 x_1,x_1)$ such that
$U(x_1-x_2)\ge\tht_2-\del_1$. The above and Lemma \ref{L.2.1}(i)
show that $\limsup_{t\to\infty} T(t,x)\le\tht_2-\del_1$ uniformly in
$x\ge \tfrac 12 x_1$, and so
\begin{equation} \lb{3.5}
\limsup_{t\to\infty} T(t,x)\le U(x-x_2)
\end{equation}
uniformly in $x\in[x_2,x_1]$. We will show that \eqref{3.5} also
holds uniformly in $x>x_1$, which will yield $x_1\le x_2$ by the
definition of $x_1$. This will be a contradiction and hence
necessarily $x_1=0$.

Let $s(t)$ be smooth, decreasing with all derivatives bounded, such
that $s(0)=\tht_2$ and $\lim_{t\to\infty} s(t)=U(x_1-x_2)$. Let
$S(t,x)$ satisfy \eqref{1.1} for $x>x_1$ with $S(0,x)=U(x-x_1)$ and
$S(t,x_1)=s(t)$. As above, one proves that this time $S$ is
non-increasing in $t$,
\begin{equation} \lb{3.6}
S(t,x)\in[U(x-x_2),U(x-x_1)],
\end{equation}
and $\til S(x)\equiv \lim_{t\to\infty} S(t,x)$ satisfies \eqref{3.3}
and $\til S(x)\to 0$ as $x\to\infty$ (by \eqref{3.6}). Moreover,
$S\to \til S$ uniformly on compacts, which together with \eqref{3.6}
and $U(x)\to 0$ as $x\to\infty$ shows that $S\to\til S$ uniformly on
$\bbR$. Since $\til S(x_1)=U(x_1-x_2)$ and $\til
S(\infty)=U(\infty)=0$, a formula similar to \eqref{3.4}, with the
integral from 0 to $U(x_1-x_2)$, gives $\til S'(x_1)=U'(x_1-x_2)$.
It follows that $\til S(x)=U(x-x_2)$.

Now pick any $\eps>0$ and
choose $t_0$ such that
\begin{equation} \lb{3.7}
S(t,x)-U(x-x_2)<\tfrac \eps 2
\end{equation}
for $t\ge t_0$ and $x\ge x_1$. Then pick $t_1$ so that
\[
T(t,x)-U(x-x_1)<\eps_0\equiv\tfrac \eps 2 e^{-ct_0}
\]
for $t\ge t_1$ and $x\ge x_1$ (with $c$ the Lipschitz constant for
$f$). This is possible by the definition of $x_1$.

For any $t_2>t_1$ and $x>x_1$ we have
\[
T(t_2,x)-S(0,x)=T(t_2,x)-U(x-x_1)<\eps_0,
\]
and for $t>t_2$ we have
\[
T(t,x_1)-S(t-t_2,x_1) = T(t,x_1)-s(t) \le T(t,x_1)-U(x_1-x_2)<\eps_0
\]
by \eqref{3.5} if $t_2$ is large enough. Hence if we let
$R(t,x)\equiv S(t,x)+e^{ct}\eps_0$, then $T(t_2,x)<R(0,x)$ for
$x>x_1$, $T(t,x_1)<R(t-t_2,x_1)$ for $t>t_2$, and
\[
R_t = S_t + ce^{ct}\eps_0 = S_{xx} + f(S) + ce^{ct}\eps_0 \ge R_{xx}
+ f(R).
\]
So $R$ is a supersolution of \eqref{1.1}, and by the comparison
principle $T(t,x)\le R(t-t_2,x)$ for $t>t_2$ and $x>x_1$. In
particular,
\[
T(t_2+t_0,x)\le R(t_0,x) = S(t_0,x)+\tfrac \eps 2 < U(x-x_2)+\eps
\]
for $x>x_1$ by \eqref{3.7}.
Since this holds for any large enough $t_2$, we have
$T(t,x)<U(x-x_2)+\eps$ for all large $t$ and $x>x_1$. As $\eps>0$
was arbitrary, this gives \eqref{3.5} uniformly in $x>x_1$. Hence
$x_1\le x_2<x_1$, a contradiction. Therefore we must have $x_1=0$
and the proof is finished.

\end{document}